\documentclass[lettersize,journal]{IEEEtran}
\usepackage{amsmath,amsfonts}
\usepackage{algorithm}
\usepackage{algorithmic}
\usepackage{array}
\usepackage{textcomp}
\usepackage{stfloats}
\usepackage{url}
\usepackage{verbatim}
\usepackage{cite}
\usepackage{tabularx}
\usepackage{amssymb}
\usepackage[font=footnotesize]{subcaption}
\usepackage[dvipdfmx]{graphicx}
\usepackage{siunitx}
\graphicspath{{./figure}}

\begin{document}

\title{Data-driven topology design for conductor layout problem of electromagnetic interference filter}

\author{Duanyutian Zhou,~Katsuya Nomura,~and Shintaro Yamasaki
\IEEEcompsocitemizethanks{
\IEEEcompsocthanksitem \textsf{\copyright 2025 IEEE. Personal use of this material is permitted. Permission from IEEE must be obtained for all other uses, in any current or future media, including reprinting/republishing this material for advertising or promotional purposes, creating new collective works, for resale or redistribution to servers or lists, or reuse of any copyrighted component of this work in other works.
}
\IEEEcompsocthanksitem Duanyutian Zhou and Shintaro Yamasaki are with the Graduate School of Information, Production and Systems, Waseda University, 2-7 Hibikino, Wakamatsu, Kitakyushu, Fukuoka 808-0135. \protect E-mail: \{zhouduanyt2000@ruri, s\underline{ }yamasaki@\}waseda.jp
\IEEEcompsocthanksitem Katsuya Nomura is with the School of Engineering, Kwansei Gakuin University, 1 Gakuen uegahara, Sanda, Hyogo 669-1330. \protect E-mail: \ katsuya.nomura@kwansei.ac.jp
}
}
\maketitle

\begin{abstract}
Electromagnetic interference (EMI) filters are used to reduce electromagnetic noise.
It is well known that the performance of an EMI filter in reducing electromagnetic noise largely depends on its conductor layout.
Therefore, if a conductor layout optimization method with a high degree of freedom is realized, a drastic performance improvement is expected.
Although there are a few design methods based on topology optimization for this purpose, these methods have some difficulties originating from topology optimization.
In this paper, we therefore propose a conductor layout design method for EMI filters on the basis of data-driven topology design (DDTD), which is a high degree of freedom structural design methodology incorporating a deep generative model and data-driven approach.
DDTD was proposed to overcome the intrinsic difficulties of topology optimization, and we consider it suitable for the conductor layout design problem of EMI filters.
One significant challenge in applying DDTD to the conductor layout design problem is maintaining the topology of the circuit diagram during the solution search.
For this purpose, we propose a simple yet efficient constraint.
We further provide numerical examples to confirm the usefulness of the proposed method.
\end{abstract}

\begin{IEEEkeywords}
electromagnetic interference filter, topology optimization, deep generative model, data-driven approach
\end{IEEEkeywords}

\section{Introduction}
\IEEEPARstart{I}{n} recent years, noise levels in power electronics circuits have been increasing due to the higher frequencies and faster switching of power devices.
International standards for power conversion circuits have been established by the International Electrotechnical Commission (IEC) and the International Special Committee on Radio Interference (CISPR), covering a range from 150 kHz to 30 MHz for conducted noise~\cite{2024_CISPR11}.
There are also standards from 9 kHz to 150 kHz in CISPR 14 and 15 for specific applications~\cite{2020_CISPR14,2018_CISPR15}.
Consequently, there is a growing need for high-performance electromagnetic interference (EMI) filters, making it essential to establish effective filter design methods.

In EMI filters, even those with the same circuit topology and components can exhibit different characteristics depending on the conductor layout~\cite{ShuoWang2004parasitic,wang2009parasitic}.
This phenomenon is attributed to parasitic elements; specifically, conductor patterns generate equivalent series inductance (ESL), and adjacent conductor patterns form parasitic capacitance~\cite{liu2019capacitive}.
Additionally, magnetic coupling between loops can significantly degrade filter performance in some cases~\cite{zeeff2003analysis,chen2008modeling,Murata2017feedthroughcapacitors}.
The dominant parasitic element varies depending on the circuit topology and parameters, making it challenging to estimate their overall impact.
The recent demand for miniaturization and higher density in filters has increased the impact of parasitic effects, further complicating the conductor design of filters.
Designing a filter that addresses multiple frequency characteristics while accounting for the effects of complex parasitic elements is a significant challenge, highlighting the need for computer-aided design methods.

Topology optimization is one of the most powerful and effective methodologies to optimize material layout.
It was originally proposed in the field of solid mechanics~\cite{bendsoe1988generating}.
The basic concept of topology optimization is to replace the original structural (or layout) design problem with a material distribution problem in a given fixed design domain.
Usually, the material distribution is represented with nodal or element-wise densities, which indicate material existence with a value from 0 to 1, on the finite element mesh of the discretized design domain.
These nodal or element-wise densities are regarded as the design variables, and they are updated towards the optimal material distribution by using mathematical programming with sensitivity information.

Due to the ability of topology optimization to obtain excellent material distributions with a high degree of design freedom, it has recently also become attractive in the field of electromagnetics, for example, in antenna design problems~\cite{Nomura2007ToEMI,zhou2010level}, electromagnetic waveguide design problems~\cite{yamasaki2011level}, and in-vehicle reactor design problems~\cite{Yamasaki2017Grayscale}.
Topology optimization has also been applied to conductor layout design problems of EMI filters~\cite{Nomura2019TOEMI, Nomura2022Density}.
In these works, the objective function is one of the scattering parameters, $S_{21}$, at a single target frequency.

Although optimized conductor layouts were successfully obtained in~\cite{Nomura2019TOEMI, Nomura2022Density}, an approximation method is used in~\cite{Nomura2022Density} to conduct the electromagnetic analysis for various conductor layouts using a fixed finite element mesh, and approximation errors caused by it are inherently unavoidable.
While one does not suffer from such approximation errors when using the conductor layout optimization method proposed in~\cite{Nomura2019TOEMI}, some issues have also been recognized through trials for further developing this method.
The first issue is unstable oscillations that occur when the topology of the conductor layout changes.
In their methods, the topology of the circuit diagram is preserved because changes to it cause catastrophic effects on the target EMI filter.
On the other hand, topology changes to the conductor layout itself are not prohibited.
However, when this occurs, the objective function sometimes significantly degrades, and the conductor material distribution attempts to restore the previous topology in the next iteration of mathematical programming.
By repeating these topology changes, unstable oscillations occur.
Once the oscillation occurs, the optimal solution search fails.

The second issue is that their methods are not suitable for minimizing $S_{21}$ at multiple target frequencies.
This issue is closely related to the first issue because the numerical instability caused by the first issue is amplified when considering two or more target frequencies.
The essential reason for these issues is that their methods search for the optimal solution using sensitivity information.
In other words, the conductor layout design problem of the EMI filter is essentially strongly nonlinear, and therefore sensitivity-based topology optimization is not suitable to solve it.

In order to resolve the essential drawback of topology optimization for strongly nonlinear problems, Yamasaki \textit{et al.}~\cite{Yamasaki2021Data} proposed a design methodology called data-driven topology design (DDTD), targeting multi-objective structural design problems.
The basic concept of DDTD is the combination of a data-driven approach and a deep generative model, a type of generative artificial intelligence.
By the combination, DDTD can generate high-performance structures (or layouts) in a sensitivity-free manner while maintaining a high degree of freedom for structural representation.

Not limited to DDTD, deep generative models have recently been often used in electromagnetic design problems.
For example, Lin \textit{et al.}~\cite{Lin2022inverse} introduced a conditional variational autoencorder (cVAE) in a scattering enhanced metasurface design problem.
They used a genetic algorithm (GA) to optimize a metasurface structure, and the cVAE was used to provide the initial individuals.
Fang \textit{et al.}~\cite{Fang2023periodic} proposed to use a variational autoencoder (VAE) in a reconfigurable periodic structure design problem.
In their method, the VAE was used to capture features of originally high-dimensional electromagnetic property data in a low dimension.
Peng and Xu~\cite{Peng2022antenna} proposed to combine a VAE and a back propagation network for accelerating the electromagnetic analysis in a double-T monopole antenna design problem.

Let us return to DDTD.
Although application studies of DDTD have been limited to structural mechanics~\cite{Yamasaki2021Data,Kii2024latent} and fluid dynamics~\cite{yaji2022data}, DDTD is considered suitable for the conductor layout design problem of the EMI filter because of its sensitivity-free and multi-objective nature.
Therefore, in this paper, we propose a conductor layout design method based on DDTD for minimizing $S_{21}$ at multiple target frequencies.

One significant challenge of this paper is how to maintain the topology of the circuit diagram under the framework of DDTD because this is a specific issue in the conductor layout design of the EMI filter.
In previous studies based on topology optimization~\cite{Nomura2019TOEMI, Nomura2022Density}, a constraint using fictitious current and electric field was introduced to maintain the topology of the circuit diagram; however, this constraint is very complicated to implement.
Instead, we introduce a simpler and more reasonable constraint to preserve it by taking advantage of DDTD.
Because DDTD can generate high-performance conductor layouts with a high degree of design freedom, it is expected that we will obtain very complex conductor layouts as the final results.
Therefore, considering the physical meanings of the obtained conductor layouts is also an important challenge in this paper.

The rest of the paper is organized as follows.
In Section~\ref{section_design_problem}, we explain a conductor layout design problem of an EMI filter, which is our target in this paper.
In Section~\ref{section_design_method}, we explain the proposed design method to solve the problem formulated in Section~\ref{section_design_problem}.
In Section~\ref{section_numerical_examples}, we provide numerical examples to demonstrate the usefulness of the proposed method and explore the physical meanings of the obtained results.
Finally, we present the conclusion in Section~\ref{section_conclusion}.

\begin{figure}[t!]
\centering
\includegraphics[scale=1.1]{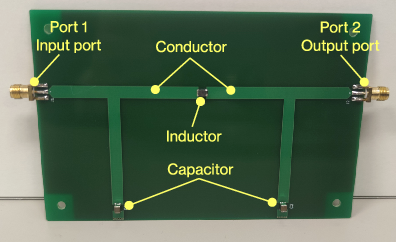}
\caption{$\pi$-type EMI filter}
\label{fig_picture}
\end{figure}

\begin{figure}[t!]
\centering
\includegraphics[scale=1.1]{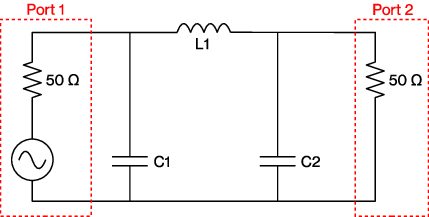}
\caption{Circuit diagram of $\pi$-type EMI filter}
\label{fig_circuit_diagram}
\end{figure}

\section{Conductor layout design problem}
\label{section_design_problem}

\subsection{$\pi$-type EMI filter}
EMI filters block adverse interferences and allow a steady flow of power to protect electronic devices and systems.
Their most important function is to reduce electromagnetic noise in a high-frequency range (typically \SI{9}{kHz} to \SI{30}{MHz}) while not affecting the passage of current in a low-frequency range (typically \SI{50}{Hz} or \SI{60}{Hz}).

Although the noise reduction performance depends on the circuit diagram design and the parameters of the electronic components, different conductor layouts lead to different parasitic effects and capacitive coupling.
Therefore, it is also important to find superior conductor layouts for improving noise reduction performance.

There are wide variations of EMI filters.
Among them, we select the $\pi$-type EMI filter shown in Fig.\ref{fig_picture} as the target in this paper.
It has two capacitors and one inductor, as shown in this figure, and these electronic components are connected to the input and output ports using conductors.
Figure~\ref{fig_circuit_diagram} shows its circuit diagram.

For the $\pi$-type EMI filter shown in Figs.~\ref{fig_picture} and \ref{fig_circuit_diagram}, the scattering parameters (so-called S-parameters) can be defined to describe the electrical behavior of the electrical network.
Among them, $S_{21}$ represents the transmission coefficient.
Therefore, we evaluate the noise reduction performance using $S_{21}$ at target frequencies.

\begin{figure}[!t]
\centering
\includegraphics[scale=1.1]{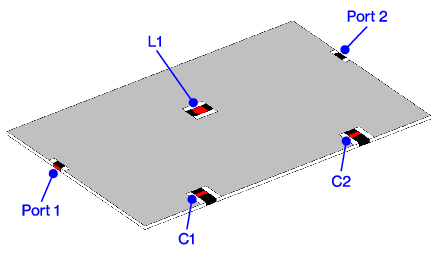}
\caption{Computation model of $\pi$-type EMI filter}
\label{fig_dd}
\end{figure}

\subsection{Numerical analysis}
Figure~\ref{fig_dd} shows the computation model.
Here, the red domains represent ports 1 and 2 (input and output ports, respectively) and the electronic components C1, C2, and L1.
The conductor material always exists in the black domains, and it does not exist in the white domains.
The back side of the printed circuit board (PCB) is filled with the conductor material.
The gray domain is the design domain, in which we consider distributions of the conductor material.
In this model, the conductor material is represented as 2D sheets using the impedance boundary condition to account for the skin effect.
The scattering boundary condition is imposed on the outer boundaries of the air domain surrounding the PCB.

Similar to~\cite{Nomura2019TOEMI}, we discretize the design domain with a numerical mesh and assign the nodal densities $\boldsymbol{\rho}$ on the mesh.
Each component of $\boldsymbol{\rho}$ takes a value from $0$ to $1$, where $0$ means there is no conductor material and $1$ means the conductor material exists.
We construct the density field $\rho$, which represents the conductor material distribution, from $\boldsymbol{\rho}$ on the basis of the shape function of the finite element method.
Then, we make a body-fitted mesh along the iso-contour of $\rho = 0.5$ for the electromagnetic analysis.
By using the body-fitted mesh, we can handle the impedance boundary condition of the conductor material accurately.
That is, we use a fixed mesh to represent various conductor material distributions and make a body-fitted mesh for each conductor material distribution when conducting the electromagnetic analysis.

We compute the electromagnetic behavior by solving the governing equations using the finite element method, as introduced in~\cite{Nomura2019TOEMI}.
We also compute $S_{21}$ at target frequencies as explained in that study.

\subsection{Formulation}
\label{section_formulation}
Although we can consider various conductor material distributions in the design domain, the topology of the circuit diagram in Fig.~\ref{fig_circuit_diagram} should be preserved.
This is an issue specific to the conductor layout design problem of the EMI filter.
In order to resolve this issue, we focus on $S_{21}$ in a low-frequency range.
As described above, $S_{21}$ represents the transmission coefficient from port~1 to port~2.
Therefore, if the disconnection of the conductor material occurs in the path from port~1 to L1, or from L1 to port~2, $S_{21}$ in a low-frequency range becomes low.
On the other hand, if the disconnection of the conductor material occurs between port~1 and C1, or port~2 and C2, $S_{21}$ in a high-frequency range becomes high.
Furthermore, if an electrical short bypassing L1 occurs, $S_{21}$ in a high-frequency range becomes high.
Therefore, we can preserve the topology of the circuit diagram by minimizing $S_{21}$ in a high-frequency range while keeping $S_{21}$ sufficiently high in a low-frequency range.

On the basis of the above discussion, we formulate the optimization problem to be solved as follows:
\begin{equation}
\label{eq_formulation}
\begin{array}{ll}
\underset{\boldsymbol{\rho}}{\text{Minimize}} & \; \left[ J_{1}(\boldsymbol{\rho}), J_{2}(\boldsymbol{\rho}) \right] \\ \\
\text{Subject to} & \; G(\boldsymbol{\rho}) \geq \bar{G}, \\ \\
& \;  0 \leq \rho_i \leq 1  \quad \text{for} \; i=1, 2, \dots, N_{\text{dsg}},
\end{array}
\end{equation}
where
\begin{equation}
\label{eq_obj_cns}
\begin{array}{l}
J_{1}(\boldsymbol{\rho}) := 20 \log_{10} \left| {S_{21}}^{(f_{1})} \right| , \\ \\
J_{2}(\boldsymbol{\rho}) := 20 \log_{10} \left| {S_{21}}^{(f_{2})} \right| , \\ \\
G(\boldsymbol{\rho}) := 20 \log_{10} \left| {S_{21}}^{(f_{3})} \right| ,
\end{array}
\end{equation}
and ${S_{21}}^{(*)}$ represents $S_{21}$ at a frequency.
Frequencies $f_{1}$ and $f_{2}$ are selected from a high-frequency range, and $f_{3}$ is selected from a low-frequency range.
As described in (\ref{eq_obj_cns}), we evaluate the decibel values of $S_{21}$ at given target frequencies as the objective and constraint functions $J_{1}(\boldsymbol{\rho})$, $J_{2}(\boldsymbol{\rho})$, and $G(\boldsymbol{\rho})$.
$\bar{G}$ is the allowable lower bound of $G(\boldsymbol{\rho})$.
In this formulation, we consider $\boldsymbol{\rho}$ as the design variables, where $\rho_i$ is the $i$-th component of $\boldsymbol{\rho}$, and $N_{\text{dsg}}$ is the number of the design variables.

\begin{figure*}[!t]
\centering
\includegraphics[scale=1.1]{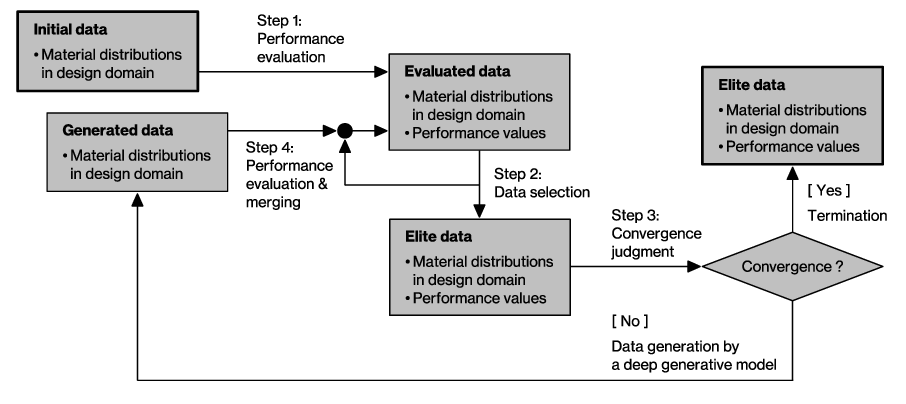}
\caption{Data process flow of the proposed method}
\label{fig_flowchart}
\end{figure*}

\section{Design method}
\label{section_design_method}
In this section, we explain the proposed design method to solve the optimization problem formulated in (\ref{eq_formulation}).

\subsection{Basic concept of DDTD}
The data process flow of the proposed method is based on DDTD.
In DDTD, structures (or layouts) are represented as material distributions in a given design domain, the same as in topology optimization, and elite (high-performance) material distributions are selected from already obtained material distribution data on the basis of the non-dominated rank\cite{deb2002fast} in the multi-objective function space.
Then, the deep generative model learns features of these elite material distributions and generates new material distributions that are diverse but inherit their features.
Some of the generated material distributions are expected to have higher performance than the current elite material distributions; therefore, the elite material distributions are updated by selecting from the merged data of the current elite and newly generated material distributions.
After that, the updated elite material distributions are used to train the deep generative model again.
By repeating these processes, the total performance of the elite material distributions is drastically improved.


\subsection{Data process flow}
\label{section_data_process_flow}
The data process flow of the proposed method is shown in Fig.~\ref{fig_flowchart}.
It starts by providing initial conductor material distributions in some way.
In this paper, we introduce a parametric model and provide initial conductor material distributions by randomly changing parameters included in the parametric model.
The initial conductor material distributions are treated as the initial data in the form of $\boldsymbol{\rho}$.
The details of the parametric model are explained in Section~\ref{section_parametric_model}.

In step~1, we evaluate the $S_{21}$ values of each conductor material distribution in the initial data using finite element analysis.
The conductor material distributions with their performance values are treated as the evaluated data.

In step~2, we calculate the non-dominated rank based on the two objective functions in (\ref{eq_formulation}), for the conductor material distributions satisfying the constraint $G(\boldsymbol{\rho}) \geq \bar{G}$.
We ignore all the conductor material distributions in the evaluated data that violate the constraint.
After that, we select the rank-one conductor material distributions as the elite data.
We further copy the elite data and store the copied data for merging in step~4.

In step~3, we check the convergence.
If the convergence criterion is satisfied, we output the elite data as the final result.
Otherwise, we train a deep generative model using the conductor material distributions of the elite data.
Here, we directly use the conductor material distributions in the form of $\boldsymbol{\rho}$ as the training data, in contrast to the previous work~\cite{Yamasaki2021Data}, where the material distributions are mapped into a regular reference domain before the training.
Through the training, features of the training data are learned in a low-dimensional latent space of the deep generative model.
Therefore, we generate new conductor material distributions in the form of $\boldsymbol{\rho}$, which are diverse but inherit features of the training data, by random sampling in the latent space.
The generated conductor material distributions are treated as the generated data.

In step~4, we evaluate the $S_{21}$ values of each conductor material distribution in the generated data, similar to step~1, and the generated data with the performance values are merged into the elite data copied in step~2.
We update the evaluated data with the merged data and return to step~2.

We iterate these processes until the convergence criterion is satisfied.

\begin{figure}[!t]
\centering
\includegraphics[scale=1.1]{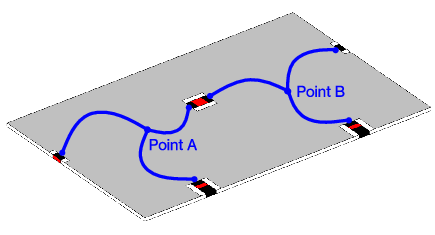}
\caption{Parametric model to provide initial conductor material distributions}
\label{fig_parametric_model}
\end{figure}

\begin{figure}[!t]
\centering
\includegraphics[scale=1.1]{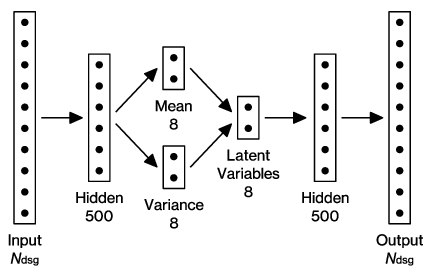}
\caption{Network architecture of VAE}
\label{fig_NN}
\end{figure}

\subsection{Parametric model for initial material distributions}
\label{section_parametric_model}
Here, we explain the parametric model to provide initial conductor material distributions.
This model represents conductor material distributions satisfying the topology of the circuit diagram by using six 2nd-order B\'{e}zier curves, as shown in Fig.~\ref{fig_parametric_model}.
We provide initial conductor material distributions by randomly changing the coordinates of points~A and B, and those of the control points of the six B\'{e}zier curves.
Note that this parametric model is used only for providing initial conductor material distributions, and the above parameters are not the design variables.
The design variables in the proposed method are the nodal densities, as explained in Section~\ref{section_formulation}.

\begin{figure}[!t]
\centering
\includegraphics[scale=0.6]{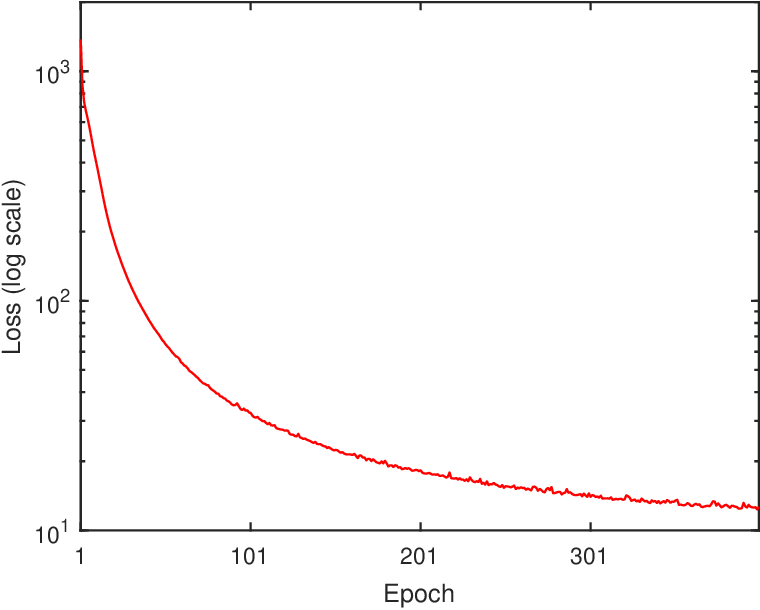}
\caption{Learning history of VAE at iteration~0 of example~1.}
\label{fig_vae_loss}
\end{figure}

\subsection{Implementation details}
In this section, we explain some implementation details.

For the generative model, we use a VAE whose network architecture is shown in Fig.~\ref{fig_NN}.
It is a simple multilayer perceptron-type architecture.
The size of the input and output layers is the number of the design variables $N_{\text{dsg}}$ in the numerical examples of Section~\ref{section_numerical_examples}.
The input layer is fully connected to the hidden layer having 500 neurons.
After the ReLU activation, these are fully connected to the 8-dimensional latent space.
It is also fully connected to the hidden layer having 500 neurons.
After the ReLU activation, these are fully connected to the output layer.
Finally, the data are output after the Sigmoid activation.

The loss function of the VAE, $L$, is given as follows:
\begin{equation}
\label{eq_loss_function}
L := L_{\text{rcn}} + \beta \cdot L_{\text{KL}},
\end{equation}
where $L_{\text{rcn}}$ is the reconstruction loss, $L_{\text{KL}}$ is the Kullback-Leibler divergence, and $\beta$ is the weight for $L_{\text{KL}}$.
Here, the reconstruction loss is measured by the mean squared error.
That is, we use the $\beta$-VAE in practice.
In this paper, we set $\beta$ to $0.5$.

The network architecture shown in Fig.~\ref{fig_NN} is very simple and almost the same as those used in design problems of structural mechanics~\cite{Yamasaki2021Data} and fluid dynamics~\cite{yaji2022data}.
As shown in Section~\ref{section_numerical_examples}, this network architecture works well without any special modifications for application to electromagnetic problems.

As described in Section~\ref{section_data_process_flow}, we train the VAE at every iteration.
This means that we have to train the VAE many times, and this is a feature of DDTD.
Therefore, we provide a constant number of training data throughout all training processes.
More precisely, we set the upper limit of the number of elite data to $400$ at step~2 of Fig.~\ref{fig_flowchart} and use the elite data as the training data.
If the number of the elite data is less than $400$, we perform data augmentation.
On the other hand, we do not use validation data.
This is because that the number of elite data is often small, particularly in early iterations, and the loss function value of the validation data cannot be used as a reliable indicator in such situations.
Instead, we train the VAE only with the training data and finish the training after $400$ epochs.
This is a VAE training manner established in~\cite{Yamasaki2021Data} and works well in subsequent studies~\cite{yaji2022data,Kii2024latent}.
For the other parameters, we set the learning rate to $1.0 \times 10^{-4}$ and the mini-batch size to $20$.
With the above settings, we train the VAE at every iteration.
As a reference, we show the learning history of the VAE at iteration~0 of example~1, which is provided in Section~\ref{section_numerical_examples}, in Fig.~\ref{fig_vae_loss}.

The trained VAE is used to generate new data.
This is done by sampling in the latent space of the trained VAE.
The sampling manner is the same as~\cite{Yamasaki2021Data}.
That is, we conduct random sampling in the range of $\left[-4, 4\right]$ for each latent variable and generate $400$ new data at every iteration.

For conductor material distributions, it is preferable that $\rho_i$ takes $0$ or $1$ except in the transition zone with a fixed width.
However, it is not ensured for conductor material distributions generated by the VAE.
Therefore, we give them this characteristic by using the normalization method proposed in~\cite{Yamasaki2021Data}.

\begin{figure*}[!t]
\centering
\includegraphics[scale=1.1]{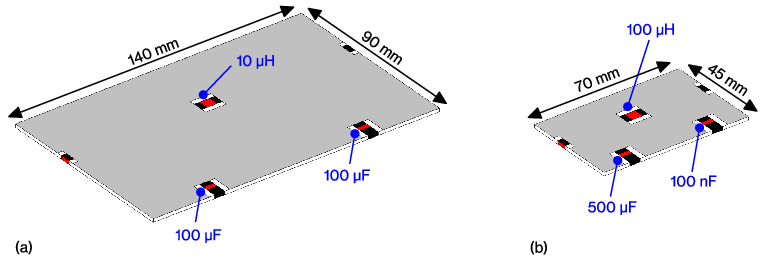}
\caption{EMI filter geometries: (a) example~1 and (b) example~2}
\label{fig_settings}
\end{figure*}

\begin{figure*}[!t]
\centering
\includegraphics[scale=1.0]{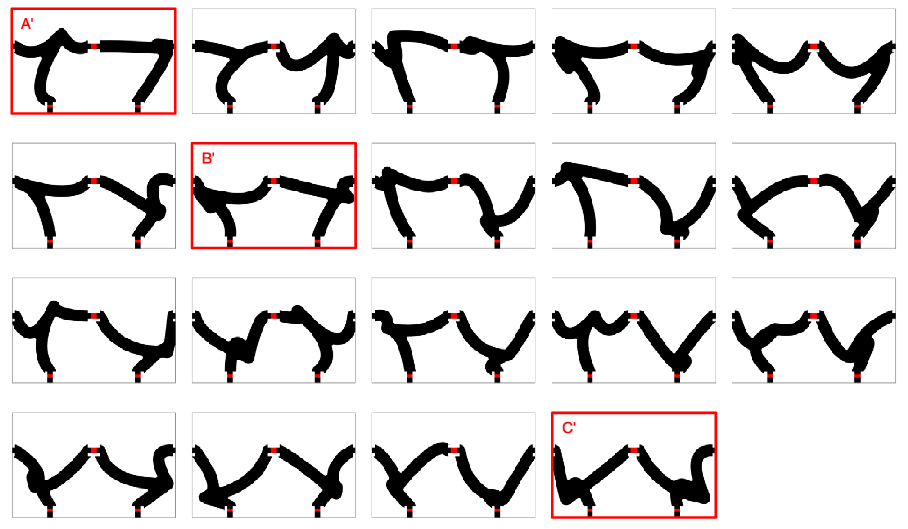}
\caption{Elite conductor layouts at iteration~0 (example~1)}
\label{fig_initial_layouts}
\end{figure*}

\begin{figure}[!t]
\centering
\includegraphics[scale=1.1]{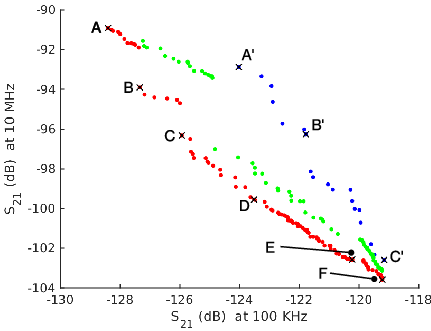}
\caption{Performance of elite conductor layouts (example~1): iteration~0 (blue), iteration~10 (green) and iteration~70 (red)}
\label{fig_result}
\end{figure}

\section{Numerical examples}
\label{section_numerical_examples}
In this section, we provide two numerical examples to demonstrate the usefulness of the proposed method.
Figure~\ref{fig_settings} shows the EMI filter geometries for these numerical examples.
In example~1, the size of the PCB is $\SI{140}{mm} \times \SI{90}{mm} \times \SI{1.6}{mm}$ (long side length, short side length, and thickness, respectively), and this PCB is surrounded by an air domain with a thickness of \SI{30}{mm}.
The constants of C1, C2 and L1 are \SI{100}{\mu F}, \SI{100}{\mu F} and \SI{10}{\mu H}, respectively.
The entire analysis domain, including the air domain, is discretized with $367,018$ finite elements, and as a result, the conductor material distribution is represented with $9,511$ nodal densities.

On the other hand, in example~2, the size of the PCB is a quarter of that of example~1; that is, $\SI{70}{mm} \times \SI{45}{mm} \times \SI{1.6}{mm}$ (long side length, short side length, and thickness, respectively), and this PCB is surrounded by an air domain with a thickness of \SI{30}{mm}.
The constants of C1, C2 and L1 are \SI{500}{\mu F}, \SI{100}{nF} and \SI{100}{\mu H}, respectively.
The entire analysis domain, including the air domain, is discretized with $329,963$ finite elements, and as a result, the conductor material distribution is represented with $8,742$ nodal densities.

We provide example~1 to confirm the usefulness of DDTD for the case in which conductive noise is dominant.
On the other hand, example~2 is for the case in which the size of the PCB becomes smaller, and inductive noise is dominant.
The details of these dominant noises are provided in the appendix.

Here, we show the common settings of examples~1 and 2.
For the material properties, the relative permittivity of the PCB is \SI{4.5}{} and that of the air is \SI{1}{}.
The relative permeability and electric conductivity are \SI{1}{} and $1 \times 10^{-8}$ {S/m} in the whole analysis domain, respectively.

In these numerical examples, we set the target frequencies $f_1$ and $f_2$ in (\ref{eq_obj_cns}) as \SI{100}{kHz} and \SI{10}{MHz}, respectively.
For the constraint to keep the topology of the circuit diagram, we found that the disconnection of the conductor material between ports~1 and 2 can be avoided by setting $f_3$ to \SI{1}{kHz} and $\bar{G}$ to \SI{-35}{dB} (example~1) or \SI{-40}{dB} (example~2).

For the convergence criterion of DDTD, we finish the computation of DDTD after 70 iterations.
We determined these parameters as a result of preliminary studies.

\subsection{Example~1}
In this section, we first show the results of DDTD for the EMI filter shown in Fig.~\ref{fig_settings}(a).
After that, we show what type of conductor material distribution, i.e., conductor layout, is obtained if we don't use the constraint to keep the topology of the circuit diagram.
Finally, we consider the reason why the obtained conductor layouts, which are very complex, have high performance from a physics viewpoint.

\begin{figure*}[!t]
\centering
\includegraphics[scale=1.0]{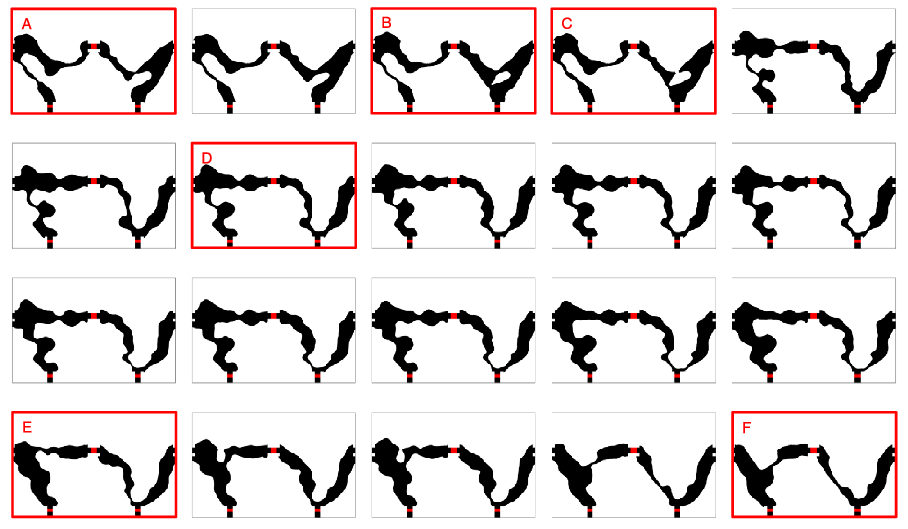}
\caption{Elite conductor layouts at iteration~70 (example~1)}
\label{fig_optimal_layouts}
\end{figure*}

\begin{figure}[!t]
\centering
\includegraphics[scale=1.0]{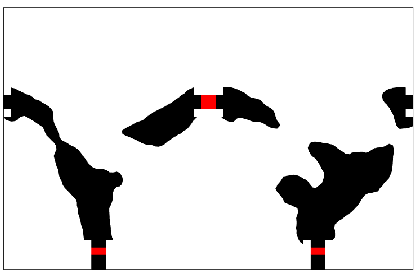}
\caption{Elite conductor layout at iteration~1 when deactivating constraint to keep the topology of the circuit diagram}
\label{fig_layouts_without_disconnection}
\end{figure}

\begin{figure}[!t]
\centering
\includegraphics[scale=0.55]{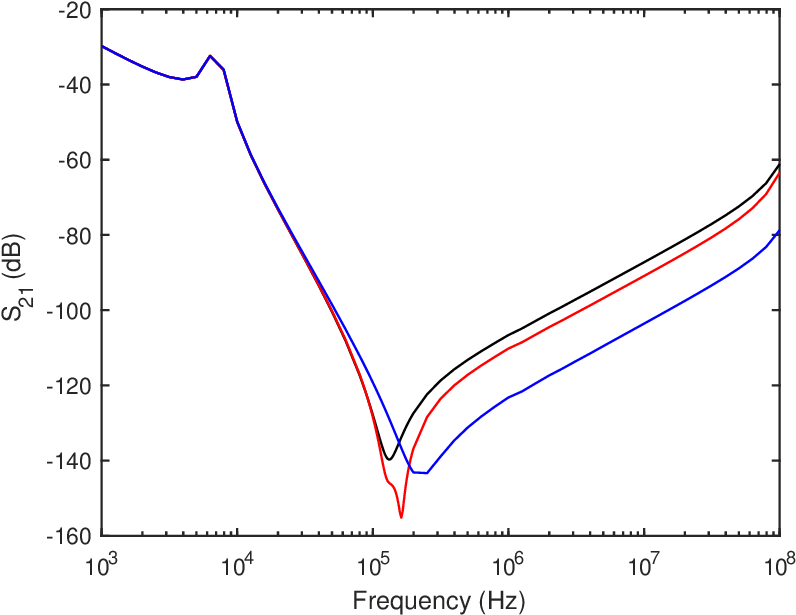}
\caption{Frequency response of $S_{21}$ (example~1): conductor layout A (red), conductor layout F (blue) and reference conductor layout (black)}
\label{fig_graph1}
\end{figure}

\begin{figure}[!t]
\centering
\includegraphics[scale=1.0]{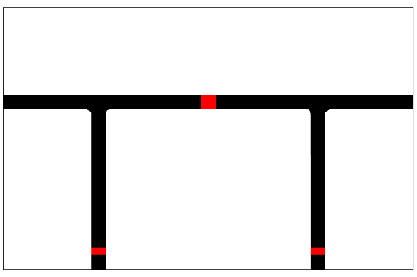}
\caption{Simple conductor layout as a reference (example~1)}
\label{fig_reference}
\end{figure}

\begin{figure}[!t]
\centering
\includegraphics[scale=1.0]{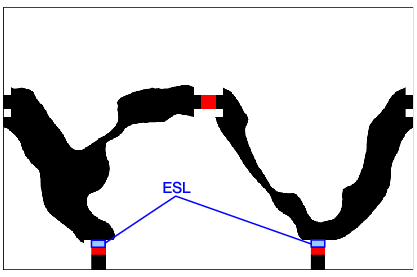}
\caption{Inductor insertion on conductor layout F (example~1)}
\label{fig_ESL}
\end{figure}

\begin{figure}[!t]
\centering
\includegraphics[scale=0.55]{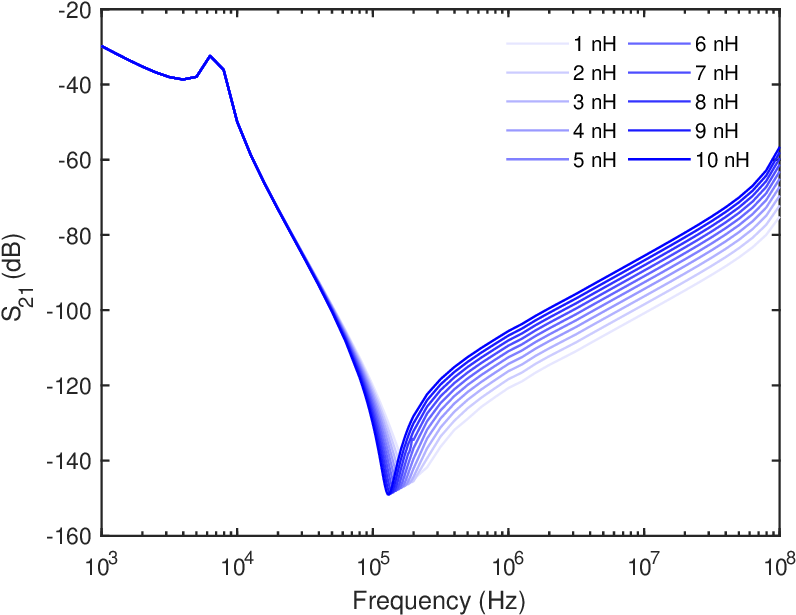}
\caption{ESL dependency on $S_{21}$}
\label{fig_graph2}
\end{figure}

\begin{figure}[!t]
\centering
\includegraphics[scale=1.0]{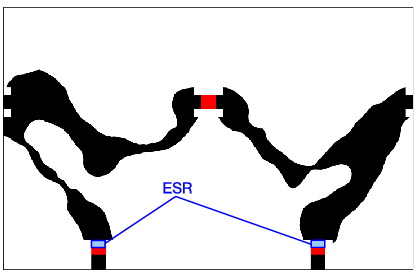}
\caption{Resistor insertion on conductor layout A (example~1)}
\label{fig_ESR}
\end{figure}

\begin{figure}[!t]
\centering
\includegraphics[scale=0.55]{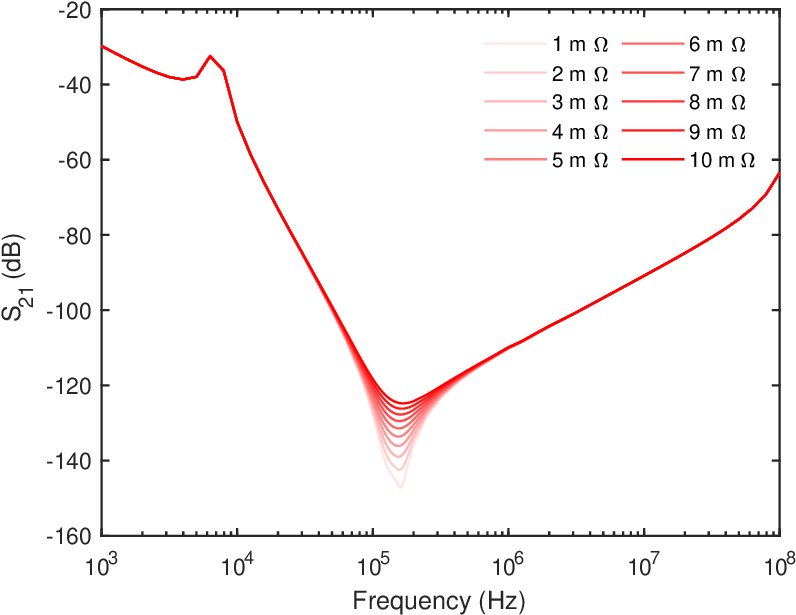}
\caption{ESR dependency on $S_{21}$}
\label{fig_graph3}
\end{figure}

\begin{figure}[!t]
\centering
\includegraphics[scale=0.55]{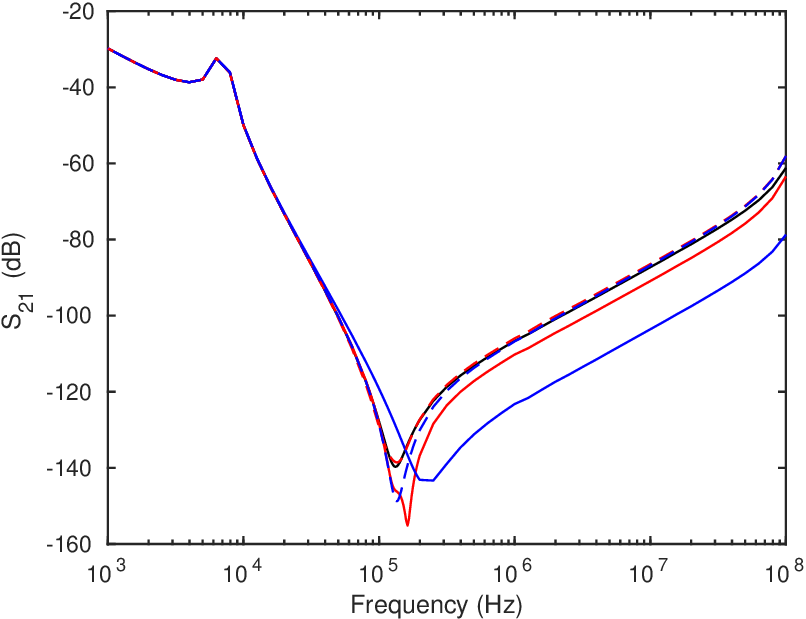}
\caption{Fitting results of ESL and ESR on $S_{21}$: reference conductor layout (black), conductor layout A (solid red), conductor layout A with ESL and ESR (dashed red), conductor layout F (solid blue), and conductor layout F with ESL (dashed blue)}
\label{fig_graph_fitting}
\end{figure}

To conduct the proposed method under these problem settings, we make $100$ initial conductor layouts as explained in Section~\ref{section_parametric_model}.
We then select elite conductor layouts on the basis of the formulation in (\ref{eq_formulation}).
The number of elite conductor layouts is $19$, and these are shown in Fig.~\ref{fig_initial_layouts}.
Starting from these conductor layouts, the performance of the elite conductor layouts are improved through the DDTD processes, as shown in Fig.~\ref{fig_result}.
The number of finally obtained elite conductor layouts is $103$, and some of them are shown in Fig.~\ref{fig_optimal_layouts}.
In Figs.~\ref{fig_initial_layouts} and \ref{fig_result}, we select three conductor layouts at iteration~0 and label them as A', B', and C'.
Similarly, in Figs.~\ref{fig_result} and \ref{fig_optimal_layouts}, we select six conductor layouts at iteration~70 and label them as A, B, C, D, E, and F.

When comparing the performance of the final elite conductor layouts (red dots in Fig.~\ref{fig_result}) to that of the initial elite conductor layouts (blue dots), the performance improvement is not significant in some combinations, e.g., F and C'.
On the other hand, we can confirm significant performance improvements in many combinations, e.g., B and A', C and A', and D and B'.
The main purpose of multi-objective optimization is to obtain many high-performance solutions in the multi-objective function space, and in general, the degree of performance improvement is not uniform.

Next, we examine the case where the constraint to keep the topology of the circuit diagram is deactivated.
Although we use the same initial conductor layouts shown in Fig.~\ref{fig_initial_layouts}, we obtain only one elite conductor layout shown in Fig.~\ref{fig_layouts_without_disconnection} at iteration~1.
Indeed, this conductor layout has excellent performance from the viewpoint of $S_{21}$ at \SI{100}{kHz} and \SI{10}{MHz}; these are \SI{-208.5}{dB} and \SI{-138.3}{dB}, respectively.
However, such an excellent conductor layout is obviously meaningless because the disconnection of the conductor material occurs between ports~1 and 2.
This numerical result indicates that the constraint we introduced is sufficiently effective to keep the topology of the circuit diagram and is essentially needed.
Here, we emphasize that our proposed constraint is simple to implement, whereas the constraint proposed in the previous work~\cite{Nomura2019TOEMI} is complicated.
We need only one electromagnetic analysis for the constraint; on the other hand, they need six electromagnetic analyses to keep the circuit diagram topology of the $\pi$-type filter in total.
This is a significant advantage of the proposed constraint.
However, the proposed constraint cannot identify disconnections when multiple parallel paths exist.
In such cases, the constraint from the previous study~\cite{Nomura2019TOEMI}, which can preserve the circuit diagram, is effective.
In other words, the two constraint methods represent a trade-off between the strictness and the ease of the implementation and can be considered complementary.

Although the elite conductor layouts at iteration~70 have superior performance as shown in Fig.~\ref{fig_result}, these layouts are complex as shown in Fig.~\ref{fig_optimal_layouts}.
Therefore, we here consider the physical reasons for the superiority of these complex layouts.
Among the $103$ layouts at iteration~70, we focus on layout~F, which has the smallest $S_{21}$ at \SI{10}{MHz}, and layout~A, which has the smallest $S_{21}$ at \SI{100}{kHz}.
For layouts~A and F, we evaluate their $S_{21}$ over a wide frequency range.
Figure~\ref{fig_graph1} shows the frequency response of $S_{21}$ for conductor layouts A and F, and also shows that of a conductor layout shown in Fig.~\ref{fig_reference} as a reference.

First, comparing the performance of conductor layout~F with that of the reference conductor layout, the characteristics are greatly improved in the high-frequency range including \SI{10}{MHz}, but conversely deteriorated at \SI{100}{kHz}.
Looking at conductor layout~F, the two branch points have both moved downward compared with the reference conductor layout, and the distance from the branch point to the capacitor has been shortened.
This is expected to reduce equivalent series inductance (ESL) and increase the bypass effect, resulting in improved performance in the high-frequency range.
To verify this hypothesis, we insert inductors at the position shown in Fig.~\ref{fig_ESL} to check the effect of ESL.
The results of the analysis with increasing ESL from \SI{1}{\nano\henry} to \SI{10}{\nano\henry} in \SI{1}{\nano\henry} increments are shown in Fig.~\ref{fig_graph2}.
As ESL increased, the characteristics in the high-frequency range, including \SI{10}{\MHz}, deteriorated, and furthermore, the characteristics at \SI{100}{\kHz} improved as the resonant frequency decreased with increasing ESL.
Therefore, both the performance change of conductor layout~F and the trade-off characteristics shown in Fig.~\ref{fig_result} were found to be due to ESL.

\begin{figure*}[!t]
\centering
\includegraphics[scale=1.0]{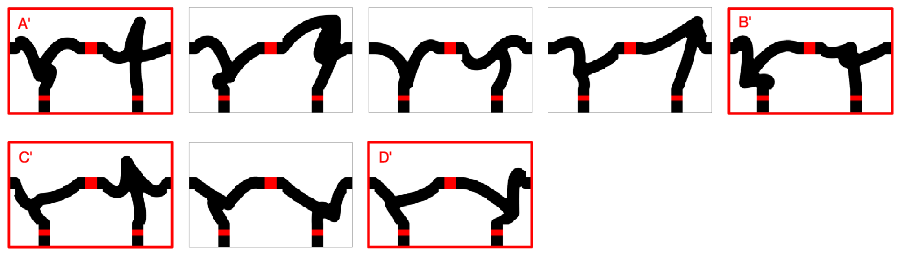}
\caption{Elite conductor layouts at iteration~0 (example~2)}
\label{fig_ex3_itr000}
\end{figure*}

\begin{figure}[!t]
\centering
\includegraphics[scale=1.1]{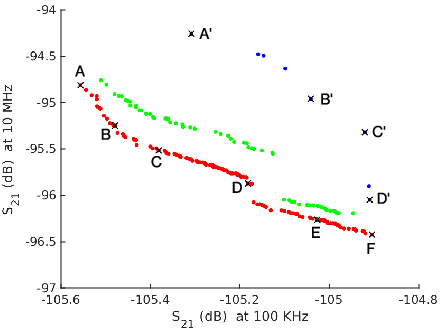}
\caption{Performance of elite conductor layouts (example~2): iteration~0 (blue), iteration~10 (green) and iteration~70 (red)}
\label{fig_ex3_result}
\end{figure}

\begin{figure*}[!t]
\centering
\includegraphics[scale=1.0]{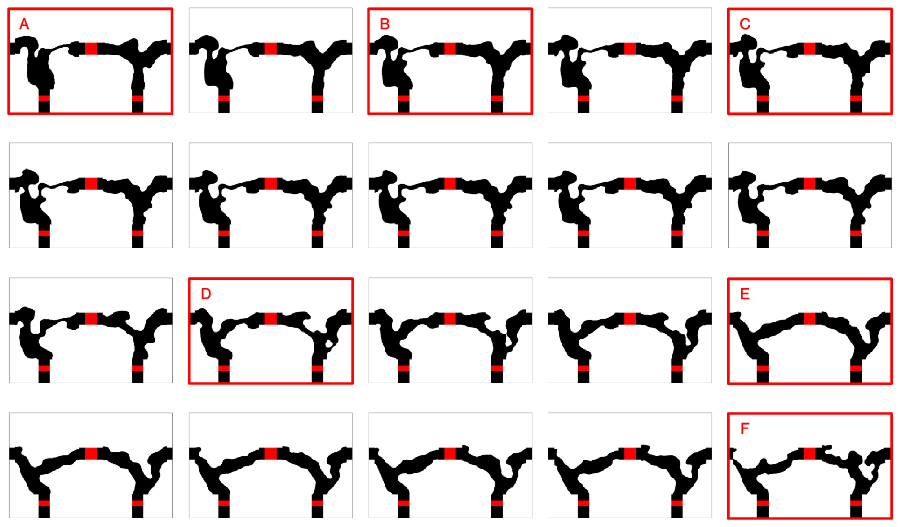}
\caption{Elite conductor layouts at iteration~70 (example~2)}
\label{fig_ex3_itr070}
\end{figure*}

\begin{figure}[!t]
\centering
\includegraphics[scale=0.6]{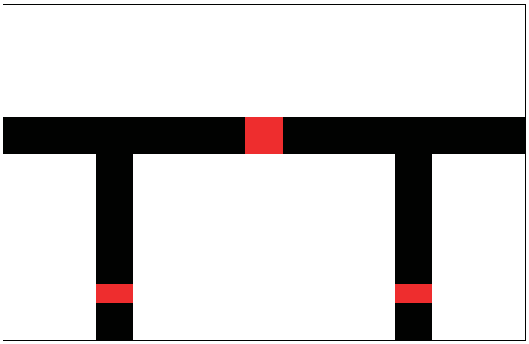}
\caption{Simple conductor layout as a reference (example 2)}
\label{fig_ex3_reference}
\end{figure}

\begin{figure}[!t]
\centering
\includegraphics[scale=0.55]{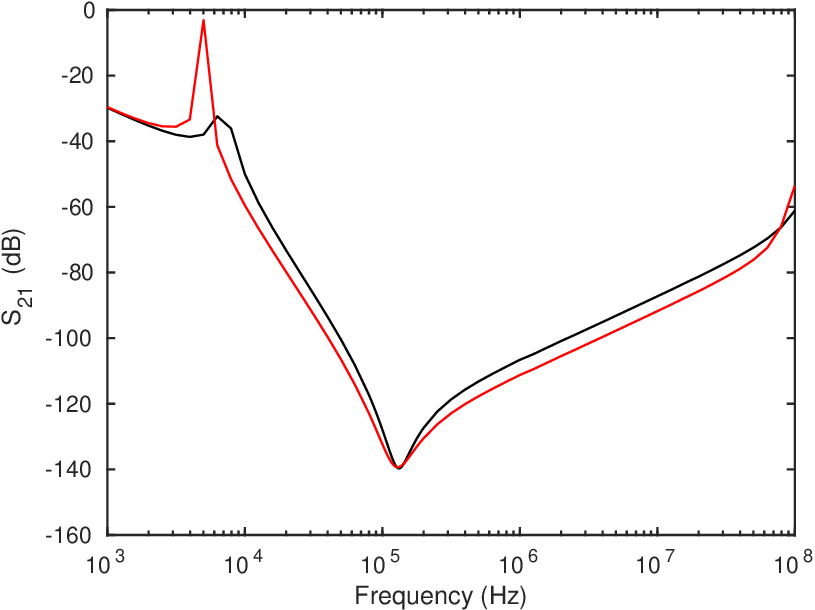}
\caption{Inductor L1 dependency on $S_{21}$ (example 1): original L1 (black) and doubled L1 (red)}
\label{fig_graph_ex1Lvar}
\end{figure}

\begin{figure}[!t]
\centering
 \includegraphics[scale=0.55]{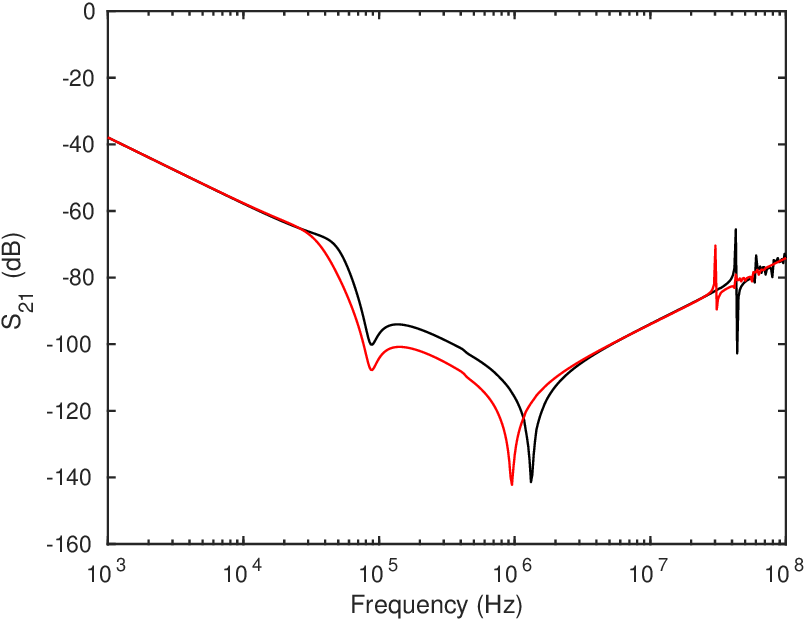}
\caption{Inductor L1 dependency on $S_{21}$ (example 2): original L1 (black) and doubled L1 (red)}
\label{fig_graph_ex2Lvar}
\end{figure}

On the other hand, in Fig.~\ref{fig_graph1}, the performance of conductor layout~A is slightly improved at high-frequencies, including \SI{10}{\MHz}, and is almost the same at \SI{100}{\kHz} compared to the reference conductor layout.
Both responses show a dip due to resonance at a frequency slightly above \SI{100}{\kHz}, and the minimum $S_{21}$ value is about \SI{-140}{\dB} for the reference conductor layout, while it is less than \SI{-150}{\dB} for conductor layout~A.
We expected this to be due to the lower equivalent series resistance (ESR) at the resonant frequency in conductor layout~A.
To verify this hypothesis, we inserted resistors at the position shown in Fig.~\ref{fig_ESR} to check the effect of ESR.
The results of the analysis with increasing ESR from \SI{1}{\mohm} to \SI{10}{\mohm} in \SI{1}{\mohm} increments are shown in Fig.~\ref{fig_graph3}.
$S_{21}$ at \SI{100}{\kHz} deteriorates as ESR increases.
Thus, it is shown that the performance improvement of conductor layout~A can be explained by the reduction of ESR in addition to the reduction of ESL.

In the above consideration, the ESL and ESR components were added in series with the capacitor within the elite conductor layout and were adjusted to fit the characteristics of the reference conductor layout.
As illustrated by the blue dashed line in Fig.~\ref{fig_graph_fitting}, inserting \SI{9}{\nano\henry} ESLs in series in the blue regions of conductor layout F (shown in Fig.~\ref{fig_ESL}) approximately reproduced the characteristics of the reference conductor layout depicted in Fig.~\ref{fig_reference}.
Similarly, as shown by the red dashed line in Fig.~\ref{fig_graph_fitting}, inserting \SI{3}{\nano\henry} ESLs and \SI{3}{\mohm} ESRs in series in the same blue regions of conductor layout A (shown in Fig.~\ref{fig_ESR}) successfully reproduced the characteristics of the reference conductor layout.
Unfortunately, the dip in the characteristics of the reference conductor layout could not be reproduced by adjusting the ESR in layout F.
However, it can be inferred that the elite conductor layouts effectively reduced the ESL and ESR corresponding to these values.

\subsection{Example~2}
In this section, we show the results of DDTD for the EMI filter shown in Fig.~\ref{fig_settings}(b).
Similar to example~1, we start by making initial conductor layouts using the parametric model explained in Section~\ref{section_parametric_model}.
Here, we make $400$ initial conductor layouts.
After that, we evaluate their performance and select $8$ elite conductor layouts under the formulation in (\ref{eq_formulation}).
These are shown in Fig.~\ref{fig_ex3_itr000}.

Figure~\ref{fig_ex3_result} shows how the performance of the elite conductor layouts is improved by DDTD.
As shown in this figure, the performance is slightly improved through $70$ iterations, and as a result, $143$ elite conductor layouts are obtained at iteration~70.
We show some of them in Fig.~\ref{fig_ex3_itr070}.

The performance improvement is smaller compared to the result shown in Fig.~\ref{fig_result} from example 1.
This is because reducing the circuit board size to one-fourth limited the impact of the conductor layout on performance.
Despite this, the proposed design method still achieved performance improvement, even in a case where the dominant noise propagation factors were different.

\section{Conclusion}
\label{section_conclusion}
In this paper, we proposed a conductor layout design method based on DDTD for the EMI filter.
DDTD is a structural design methodology with a high degree of freedom similar to topology optimization, but in contrast to it, DDTD is sensitivity-free and for multi-objective optimization problems.
By leveraging DDTD's nature, we solved a conductor layout design problem for the EMI filter, which is strongly nonlinear and multi-objective.
One important task when solving this problem is to maintain the topology of the circuit diagram during the solution search.
For this task, we proposed a simple but efficient constraint based on the $S_{21}$ value in a low-frequency range.
Including this constraint, the validity of the proposed method was confirmed using numerical examples.
As shown in the numerical examples, the obtained high-performance conductor layouts are very complex compared to a conventional conductor layout as a reference.
Therefore, we considered why the obtained conductor layouts are high-performance.
As a result, we found that the obtained conductor layouts are substantially effective in reducing ESL and ESR.

When considering the practical usage of the EMI filter, $S_{21}$ should be minimized over a given frequency range, not just at a fixed target frequency.
Therefore, as future work, we will try to improve the proposed method to minimize $S_{21}$ across an entire given frequency range.
By doing so, we will be able to avoid narrowband resonances.
It is also expected that higher-performance conductor layouts can be obtained by introducing a mutation operation, as discussed in~\cite{yaji2022data}.

Another unresolved issue is suppressing conductor layouts with narrow parts for manufacturability.
As shown in Figs.~\ref{fig_optimal_layouts} and \ref{fig_ex3_itr070}, we obtained many conductor layouts with narrow parts; however, such narrow parts may cause difficulties in manufacturing.
To address this, we have a plan to introduce a geometric constraint on the allowable minimum length into DDTD.
More specifically, we will use the signed distance function corresponding to the conductor layout to implement the geometric constraint, because its usefulness for detecting narrow parts has already been confirmed in several papers on topology optimization~\cite{2014_Guo,2015_Xia,2016_Allaire}.

We will attempt these challenges in future work.

\appendix
\section*{}
Here, we investigate noise propagation paths in the design problems of examples~1 and 2.
For the reference conductor layouts shown in Figs.~\ref{fig_reference} and \ref{fig_ex3_reference}, we evaluated the impact on $S_{21}$ when the inductance value of inductor L1, located at the center of the board, was varied.  

In example 1, doubling the inductance of L1 from \SI{10}{\micro\henry} to \SI{20}{\micro\henry}, as illustrated in Fig.~\ref{fig_graph_ex1Lvar}, resulted in a reduction of $S_{21}$ by approximately \SI{6}{\dB} in the high-frequency band, including \SI{10}{\MHz}.
This clearly indicates that the primary noise in this case is conducted through L1.

Conversely, in example 2, as shown in Fig.~\ref{fig_graph_ex2Lvar}, doubling L1 from \SI{100}{\micro\henry} to \SI{200}{\micro\henry} did not result in any significant change in $S_{21}$ at \SI{10}{\MHz}.
This suggests that the dominant noise is no longer conducted through L1 in the modified design.
It is known that the magnetic field coupling between the input and output loops is an issue in this $\pi$-type filter, particularly when the inductance of L1 is large~\cite{zeeff2003analysis,chen2008modeling,Murata2017feedthroughcapacitors}.

\bibliographystyle{IEEEtran}

\begin{thebibliography}{10}
\providecommand{\url}[1]{#1}
\csname url@rmstyle\endcsname
\providecommand{\newblock}{\relax}
\providecommand{\bibinfo}[2]{#2}
\providecommand\BIBentrySTDinterwordspacing{\spaceskip=0pt\relax}
\providecommand\BIBentryALTinterwordstretchfactor{4}
\providecommand\BIBentryALTinterwordspacing{\spaceskip=\fontdimen2\font plus
\BIBentryALTinterwordstretchfactor\fontdimen3\font minus
  \fontdimen4\font\relax}
\providecommand\BIBforeignlanguage[2]{{%
\expandafter\ifx\csname l@#1\endcsname\relax
\typeout{** WARNING: IEEEtran.bst: No hyphenation pattern has been}%
\typeout{** loaded for the language `#1'. Using the pattern for}%
\typeout{** the default language instead.}%
\else
\language=\csname l@#1\endcsname
\fi
#2}}

\bibitem{2024_CISPR11}
\emph{Industrial, scientific and medical equipment - Radio-frequency
  disturbance characteristics - Limits and methods of measurement (CISPR 11)},
  {Comit\'{e} international sp\'{e}cial des perturbations radio\'{e}lectriques
  (CISPR), Edition 7.0}, 2024.

\bibitem{2020_CISPR14}
\emph{Requirements for household appliances, electric tools and similar
  apparatus - Part 1: Emission (CISPR 14-1)}, {Comit\'{e} international
  sp\'{e}cial des perturbations radio\'{e}lectriques (CISPR), Edition 7.0},
  2020.

\bibitem{2018_CISPR15}
\emph{Limits and methods of measurement of radio disturbance characteristics of
  electrical lighting and similar equipment (CISPR 15)}, {Comit\'{e}
  international sp\'{e}cial des perturbations radio\'{e}lectriques (CISPR),
  Edition 9.0}, 2018.

\bibitem{ShuoWang2004parasitic}
S.~Wang, F.~Lee, D.~Chen, and W.~G. Odendaal, ``Effects of parasitic parameters
  on {EMI} filter performance,'' \emph{{IEEE} Transactions on Power
  Electronics}, vol.~19, no.~3, pp. 869--877, 2004.

\bibitem{wang2009parasitic}
S.~Wang, Y.~Y. Maillet, F.~Wang, R.~Lai, F.~Luo, and D.~Boroyevich, ``Parasitic
  effects of grounding paths on common-mode {EMI} filter's performance in power
  electronics systems,'' \emph{{IEEE} Transactions on Industrial Electronics},
  vol.~57, no.~9, pp. 3050--3059, 2009.

\bibitem{liu2019capacitive}
B.~Liu, R.~Ren, F.~F. Wang, D.~Costinett, and Z.~Zhang, ``Capacitive coupling
  in {EMI} filters containing {T}-shaped joint: Mechanism, effects, and
  mitigation,'' \emph{{IEEE} Transactions on Power Electronics}, vol.~35,
  no.~3, pp. 2534--2547, 2019.

\bibitem{zeeff2003analysis}
T.~M. Zeeff, T.~H. Hubing, T.~P. van Doren, and D.~Pommerenke, ``Analysis of
  simple two-capacitor low-pass filters,'' \emph{{IEEE} Transactions on
  Electromagnetic Compatibility}, vol.~45, no.~4, pp. 595--601, 2003.

\bibitem{chen2008modeling}
H.~Chen, Z.~Qian, Z.~Zeng, and C.~Wolf, ``Modeling of parasitic inductive
  couplings in a pi-shaped common mode {EMI} filter,'' \emph{{IEEE}
  Transactions on Electromagnetic Compatibility}, vol.~50, no.~1, pp. 71--79,
  2008.

\bibitem{Murata2017feedthroughcapacitors}
Y.~Murata, K.~Takahashi, T.~Kanamoto, and M.~Kubota, ``Analysis of parasitic
  couplings in {EMI} filters and coupling reduction methods,'' \emph{{IEEE}
  Transactions on Electromagnetic Compatibility}, vol.~59, no.~6, pp.
  1880--1886, 2017.

\bibitem{bendsoe1988generating}
M.~P. Bends{\o}e and N.~Kikuchi, ``Generating optimal topologies in structural
  design using a homogenization method,'' \emph{Computer Methods in Applied
  Mechanics and Engineering}, vol.~71, no.~2, pp. 197--224, 1988.

\bibitem{Nomura2007ToEMI}
T.~Nomura, K.~Sato, K.~Taguchi, T.~Kashiwa, and S.~Nishiwaki, ``Structural
  topology optimization for the design of broadband dielectric resonator
  antennas using the finite difference time domain technique,''
  \emph{International Journal for Numerical Methods in Engineering}, vol.~71,
  no.~11, pp. 1261--1296, 2007.

\bibitem{zhou2010level}
S.~Zhou, W.~Li, and Q.~Li, ``Level-set based topology optimization for
  electromagnetic dipole antenna design,'' \emph{Journal of Computational
  Physics}, vol. 229, no.~19, pp. 6915--6930, 2010.

\bibitem{yamasaki2011level}
S.~Yamasaki, T.~Nomura, A.~Kawamoto, K.~Sato, and S.~Nishiwaki, ``A level
  set-based topology optimization method targeting metallic waveguide design
  problems,'' \emph{International Journal for Numerical Methods in
  Engineering}, vol.~87, no.~9, pp. 844--868, 2011.

\bibitem{Yamasaki2017Grayscale}
S.~Yamasaki, A.~Kawamoto, A.~Saito, M.~Kuroishi, and K.~Fujita,
  ``Grayscale-free topology optimization for electromagnetic design problem of
  in-vehicle reactor,'' \emph{Structural and Multidisciplinary Optimization},
  vol.~55, no.~3, pp. 1079--1090, 2017.

\bibitem{Nomura2019TOEMI}
K.~Nomura, S.~Yamasaki, K.~Yaji, H.~Bo, A.~Takahashi, T.~Kojima, and K.~Fujita,
  ``Topology optimization of conductors in electrical circuit,''
  \emph{Structural and Multidisciplinary Optimization}, vol.~59, no.~6, pp.
  2205--2225, 2019.

\bibitem{Nomura2022Density}
K.~Nomura, ``Density-based topology optimization for conductor design of emi
  filters with improved impedance boundary condition,'' in \emph{Proceedings of
  2022 International Symposium on Electromagnetic Compatibility - {EMC}
  Europe}, 2022, pp. 377--382.

\bibitem{Yamasaki2021Data}
S.~Yamasaki, K.~Yaji, and K.~Fujita, ``Data-driven topology design using a deep
  generative model,'' \emph{Structural and Multidisciplinary Optimization},
  vol.~64, no.~3, pp. 1401--1420, 2021.

\bibitem{Lin2022inverse}
H.~Lin, J.~Hou, J.~jin, Y.~Wang, R.~Tang, X.~Shi, Y.~Tian, and W.~Xu,
  ``Machine-learning-assisted inverse design of scattering enhanced
  metasurface,'' \emph{Optics Express}, vol.~30, no.~2, pp. 3076--3088, 2022.

\bibitem{Fang2023periodic}
X.~Fang, H.~Li, and Q.~Cao, ``Design of reconfigurable periodic structures
  based on machine learning,'' \emph{{IEEE} Transactions on Microwave Theory
  and Techniques}, vol.~71, no.~8, pp. 3341--3351, 2023.

\bibitem{Peng2022antenna}
K.~Peng and F.~Xu, ``Optimization of antenna performance based on
  {VAE-BPNN-PCA},'' in \emph{Proceedings of 2022 International Conference on
  Microwave and Millimeter Wave Technology (ICMMT)}, 2022, pp. 1--3.

\bibitem{Kii2024latent}
T.~Kii, K.~Yaji, K.~Fujita, Z.~Sha, and C.~C. Seepersad, ``Latent crossover for
  data-driven multifidelity topology design,'' \emph{Journal of Mechanical
  Design}, vol. 146, pp. 051\,713--1, 2024.

\bibitem{yaji2022data}
K.~Yaji, S.~Yamasaki, and K.~Fujita, ``Data-driven multifidelity topology
  design using a deep generative model: Application to forced convection heat
  transfer problems,'' \emph{Computer Methods in Applied Mechanics and
  Engineering}, vol. 388, p. 114284, 2022.

\bibitem{deb2002fast}
K.~Deb, A.~Pratap, S.~Agarwal, and T.~Meyarivan, ``A fast and elitist
  multiobjective genetic algorithm: {NSGA-II},'' \emph{{IEEE} transactions on
  evolutionary computation}, vol.~6, no.~2, pp. 182--197, 2002.

\bibitem{2014_Guo}
X.~Guo, W.~Zhang, and W.~Zhong, ``Explicit feature control in structural
  topology optimization via level set method,'' \emph{Computer Methods in
  Applied Mechanics and Engineering}, vol. 272, pp. 354--378, 2014.

\bibitem{2015_Xia}
Q.~Xia and T.~Shi, ``Constraints of distance from boundary to skeleton: {For}
  the control of length scale in level set based structural topology
  optimization,'' \emph{Computer Methods in Applied Mechanics and Engineering},
  vol. 295, pp. 525--542, 2015.

\bibitem{2016_Allaire}
G.~Allaire, F.~Jouve, and G.~Michailidis, ``Thickness control in structural
  optimization via a level set method,'' \emph{Structural and Multidisciplinary
  Optimization}, vol.~53, no.~6, pp. 1349--1382, 2016.

\end{thebibliography}

\end{document}